\documentclass[secnumarabic,amssym,amsmath,amsfonts,nobibnotes,twoside,aip,APS]{revtex4}

\begin{document}

\title{\bf On  Agaian Matrix }
\author{P. Di\c t\u a}
\email{dita@zeus.theory.nipne.ro}

\affiliation{ Horia Hulubei National Institute of Physics and Nuclear Engineering, P.O. Box MG6, Magurele,  Romania }

\begin{abstract}
In this paper we discuss  the isolation problem of the Agaian matrix among the 6-dimensional matrices. In general a matrix isolation could depend on the kind of equivalence one uses for complex Hadamard matrices. However the present form of Agaian matrix has  such a big symmetry that both  equivalence methods lead to the same result: the matrix is isolated.
\end{abstract}

\maketitle

\section{Introduction}

One of the most interesting 6-dimensional matrices is the Agaian matrix
which appears in \cite{A}  on page 112. By using the usual equivalence of complex Hadamard matrices this matrix is an isolated one. It has the form
\begin{eqnarray}
A_1(\omega) = \left[ \begin{array}{llllll}
1&1&1&1&1&1\\*[2mm]
1&1&\omega&\omega^2&\omega^2&\omega\\*[2mm]
1&\omega&1&\omega&\omega^2&\omega^2\\*[2mm]
1&\omega^2&\omega&1&\omega&\omega^2\\*[2mm]
1&\omega^2&\omega^2&\omega&1&\omega\\*[2mm]
1&\omega&\omega^2&\omega^2&\omega&1\end{array}\right]\label{Ag2}\end{eqnarray}
where $1,\,\,\omega\,\, {\rm  and}\,\, \omega^2$ are the solutions of equation $x^3-1=0$, and $i = \sqrt{-1}$. 

Nowadays one makes use of two different equivalence methods: the standard method and the unitary equivalence.
The first one has its origin in Sylvester paper \cite {Sy} who introduced the so called standard form for real Hadamard matrices, which implies that the entries of the first row and column are equal to 1. This equivalence was also extended  to complex Hadamard matrices. It says that two Hadamard matrices $H_1$ and $H_2$ are  equivalent, written as $H_1\equiv H_2$, if there exist two diagonal unitary matrices $D_1$ and $D_2$, and permutation matrices $P_1$ and $P_2$, such that 
\begin{eqnarray}
H_1=D_1P_1H_2P_2D_2\label{Ha}\end{eqnarray}

However the complex Hadamard matrices naturally belong to the class of  normal matrices. A matrix $N$ is normal if it commutes with its adjoint $N^*$, i.e. it satisfies the relation $N\,N^*= N^*\,N$. For this class of operators the unitary equivalence takes a simple form and says that every normal matrix is similar to a diagonal matrix $D$, which means that there exists a unitary matrix $U$, such that
\begin{eqnarray} N= UDU^* \label{no}\end{eqnarray}
 see \cite{P} p. 357. For all unitary matrices the entries of the diagonal matrix $D$ are unimodular.

 Two important classes of normal operators are the unitary and self-adjoint matrices, such that for both these classes two matrices $M_1$ and $M_2$  are unitary equivalent
 iff they have the same spectrum, or equivalently, the characteristic polynomials are the same up to a multiplicative constant factor, see  \cite{P}, i.e.
\begin{eqnarray}
pol(M_1(x))= det(x\,I_n-M_1/\sqrt{n})=det(x\,I_n-M_2/\sqrt{n}=pol(M_2(x))\end{eqnarray}
where $det$ is the determinant of the corresponding matrix. In the following we shall make use of this type of equivalence. To understand better the difference between the two equivalence methods  we give a simple example of a selfadjoint Hadamard matrix
\begin{eqnarray}
M_6=\left[ \begin{array}{rrrrrr}
1&1&1&1&1&1\\
1&-1&i&i&-i&-i\\
1&-i&-1&1&-1&i\\
1&-i&1&-1&i&-1\\
1&i&-1&-i&1&-1\\
1&i&-i&-1&-1&1
\end{array}\right]\label{nor}\end{eqnarray}

The unitary equivalence implies the knowledge of its spectrum, which in this case is given by
\begin{eqnarray} Sp(M_6)=\left[-1^3,1^3\right]\label{spm}\end{eqnarray}
where power means eigenvalue multiplicity. The above result shows that $M_6$ is a complex Hadamard matrix, and in the same time a selfadoint matrix, its entries and  spectrum being unimodular and real.

If we make use  of the standard equivalence (\ref{Ha})  the following matrix 
\begin{eqnarray}
M_{61}=\left[ \begin{array}{rrrrrr}
1&1&1&1&1&1\\
1&-1&1&-1&i&-i\\
1&1&-1&i&-1&-i\\
1&-i&-1&-1&1&i\\
1&-1&-i&1&-1&i\\
1&i&i&-i&-i&-1
\end{array}\right]\label{nor1}\end{eqnarray}
 is equivalent to (\ref{nor}), but its spectrum  is given by
\begin{eqnarray} Sp(M_{61})=\left[-1^2,1^2,\frac{i -\sqrt{2}}{\sqrt{3}} ,-\frac{i+\sqrt{2}}{\sqrt{3}}\right]\label{spm1}\end{eqnarray}
Thus when  one makes use of the usual equivalence, see  relation (\ref{Ha}), the classical quantum mechanics will be in a big danger, because the $M_{61}$ eigenvalues are unimodular, but not real.

In the first case $M_{6} = M_{6}^*$, and in the second case  $M_{61} \ne M_{61}^*$.  Thus the change of rows and/or colums between themselves modifies the matrix symmetry, and accordingly its spectrum.

\section{ Agaian Matrix}

The Agaian matrix, (\ref{Ag2}),  is a Butson type matrix, see papers \cite{B1}-\cite{B2}, and its first form appeared on page 104, in book \cite{A}, Chapter 2, paragraph 5, named Generalized Hadamard Matrices.  Agaian notation is
\begin{eqnarray}
H(3,6)=\left[ \begin{array}{rrrrrr}
z&x&y&y&x&z\\
x&z&y&x&y&z\\
x&x&x&x&x&x\\
z&x&z&x&y&y\\
x&z&z&y&x&y\\
z&z&x&y&y&x\end{array}\right]\label{Ag}\end{eqnarray}
where $x=1,\,\, y=\omega,\,\,z=\omega^2$ are the third roots of unity, and it does not coincide with matrices from his note 5.1 where Agaian gave a few examples of generalized Hadamard matrices. 

After providing the above matrix (\ref{Ag}) Agaian said that each generalized Hadamard matrix can be brought to a standard form matrix, i.e. with first row and column entries equal to 1. After that he wrote:

{\em Butson (1962) proved that for prime numbers $p$ the necessary condition for existence of normalized H(p,h)-matrices is $h=p t$ where $t$ is a natural number}.

 In fact in both the papers, \cite{B1} and \cite{B2}, are given only  necessary conditions for the existence of a matrix whose entries come from  Fourier matrix, and no numerical example is provided. Im our opinion this means that the above matrix was firstly obtained by Agaian.

With the above notations  the $H(3,6)$ matrix (\ref{Ag}) has the form 

\begin{eqnarray}
A_{10}=\left[ \begin{array}{rrrrrr}
\omega^2&1&\omega&\omega &1&\omega^2 \\*[2mm]
1&\omega^2 &\omega&1&\omega&\omega^2 \\*[2mm]
1&1&1&1&1&1\\*[2mm]
\omega^2&1&\omega^2&1&\omega&\omega\\*[2mm]
1&\omega^2 &\omega^2 &\omega&1&\omega\\*[2mm]
\omega^2 &\omega^2 &1&\omega&\omega& 1\end{array}\right]\label{Ag1}\end{eqnarray}
Because we make use of unitary equivalence we compute the spectral function of $A_{10}$-matrix and its form is given by

\begin{eqnarray}
f(A_{10})= x^6+\sqrt{\frac{2}{3}}(\omega-1)x^5-\frac{1+2\omega}{2}x^4+\frac{1+2\omega}{\sqrt{6}}x^3-\frac{1+2\omega}{2}x^2+\sqrt{\frac{2}{3}}(2+\omega)x -1
\label{sp10}\end{eqnarray}
 whose spectrum cannot be explicitly found because the above equation is not a reciprocal one. However its spectrum can be calculated numerically by using Mathematica, and its spectrum is simple.

Multiplying at left and right by unitary diagonal matrices generated by the first row and/or the first column the $A_{10}$ matrix transforms into a standard form  matrix 
\begin{eqnarray}
A_{01}=\left[ \begin{array}{llllll}
1&1&1&1&1&1\\*[2mm]
1&\omega&\omega^2&\omega&1&\omega^2\\*[2mm]
1&\omega^2&\omega&\omega&\omega^2&1\\*[2mm]
1&1&\omega&\omega^2&\omega&\omega^2\\*[2mm]
1&\omega&1&\omega^2&\omega^2&\omega\\*[2mm]
1&\omega^2&\omega^2&1&\omega&\omega \end{array}\right]\label{Ag3}\end{eqnarray}
whose spectral function is simpler than (\ref{sp10}) and has the form
\begin{eqnarray}
f(A_{01})=(x^2-1)(x^4+\frac{1-\omega}{\sqrt{6}}x^3+x^2+\frac{\omega+2}{\sqrt{6}}x+1) \label{sp01}\end{eqnarray}
In this case its spectrum can be easily found. 

Perhaps the above form (\ref{Ag3}) was considered as being not sufficiently symmetric, and it was easily  transformed into matrix (\ref{Ag2}) whose spectrum does not depend on $\omega$, and is given by
\begin{eqnarray}
Sp(A_1(\omega))=\left[-1,1,\left(\frac{\sqrt{3}-i \sqrt{5}}{2 \sqrt{2}}\right)^2,\left(\frac{\sqrt{3}+i \sqrt{5}}{2 \sqrt{2}}\right)^2 \right]\label{Sp(A1)}\end{eqnarray}

Looking at matrices (\ref{Ag}) and  (\ref{Ag1}) one observes that Again used only a particular case, namely that when $x=1,\,y=\omega,\,z=\omega^2$. However the three parameters, (x,\,y,\,z), generate six permutations and this suggests that there could be at least six different matrices, if one makes use of unitary equivalence, instead of the usual one.

For example if we take $x=\omega,\,y=1,\, z=\omega^2$ in matrix (\ref{Ag}) the form of the corresponding matrix is
\begin{eqnarray}
A_{20}=\left[ \begin{array}{rrrrrr}
\omega^2&\omega&1&1&\omega &\omega^2 \\*[2mm]
\omega &\omega^2&1&\omega&1&\omega^2 \\*[2mm]
\omega&\omega&\omega& \omega&\omega&\omega\\*[2mm]
\omega^2&\omega&\omega^2&\omega&1&1\\*[2mm]
\omega &\omega^2 &\omega^2&1&\omega&1\\*[2mm]
\omega^2 &\omega^2 &\omega&1&1&\omega\end{array}\right]\label{Ag5}\end{eqnarray}

One sees that $A_{10}\neq  A_{20}$.  In fact the spectral function of (\ref{Ag5}) has the form 
\begin{eqnarray}
f(A_{20})= x^6+\sqrt{\frac{2}{3}}(1-\omega)x^5+\frac{1-\omega}{2}x^4-\frac{1+2\omega}{\sqrt{6}} x^3-\frac{\omega+2}{2}x^2-\sqrt{\frac{2}{3}}(\omega+2)x -1
\label{sp20}\end{eqnarray}

The corresponding standard form is 
\begin{eqnarray}
A_{02}=\left[ \begin{array}{llllll}
1&1&1&1&1&1\\*[2mm]
1&\omega^2&\omega&\omega^2&1&\omega\\*[2mm]
1&\omega&\omega^2&\omega^2&\omega&1\\*[2mm]
1&1&\omega^2&\omega&\omega^2&\omega\\*[2mm]
1&\omega^2&1&\omega&\omega&\omega^2\\*[2mm]
1&\omega&\omega&1&\omega^2&\omega^2 \end{array}\right]\label{Ag7}\end{eqnarray}
and its spectral function is given by
\begin{eqnarray}
f(A_{02})=(x^2-1)(x^4+\frac{\omega+2}{\sqrt{6}} x^3+x^2+\frac{1-\omega}{\sqrt{6}}x+1)
\label{sp02}\end{eqnarray}

Another form of the matrix (\ref{Ag7}) is the following

\begin{eqnarray}
A_{2}(\omega)=\left[ \begin{array}{rrrrrr}
1&1&1&1&1&1\\*[2mm]
1&1&\omega^2&\omega&\omega^2&\omega\\*[2mm]
1&\omega^2&1&\omega&\omega&\omega^2\\*[2mm]
1&\omega&\omega &1&\omega^2&\omega^2\\*[2mm]
1&\omega^2&\omega&\omega&1&\omega^2\\*[2mm]
1&\omega^2&\omega&\omega^2&\omega&1\end{array}\right]\label{Ag6}\end{eqnarray}
and it does not coincide with the $A_1$-matrix, see (\ref{Ag2}); however its spectrum is given by the same relation (\ref{Sp(A1)}).

For the choice $x = \omega,\,y = \omega^2,\,z = 1$ one finds
\begin{eqnarray}
A_{30}=\left[ \begin{array}{rrrrrr}
1&\omega&\omega^2&\omega^2 &\omega&1 \\*[2mm]
\omega &1&\omega^2&\omega&\omega^2&1 \\*[2mm]
\omega&\omega&\omega& \omega&\omega&\omega\\*[2mm]
1&\omega&1&\omega&\omega^2&\omega^2\\*[2mm]
\omega &1&1&\omega^2 &\omega&\omega^2 \\*[2mm]
1&1&\omega &\omega^2 &\omega^2&\omega\end{array}\right]\label{Ag8}\end{eqnarray}
and its spectral function is
\begin{eqnarray}
f(A_{30})= x^6-\sqrt{\frac{2}{3}}(1+2\omega)x^5+\frac{\omega-1}{2}x^4+\frac{1+2\omega}{\sqrt{6}} x^3+\frac{\omega+2}{2}x^2-\sqrt{\frac{2}{3}}(1+2\omega)x -1
\label{sp30}\end{eqnarray}

Similar to the preceding cases one gets 
\begin{eqnarray}
A_{03}=\left[ \begin{array}{rrrrrr}
1&1&1&1&1&1\\*[2mm]
1&\omega&\omega^2&\omega&1&\omega^2\\*[2mm]
1&\omega^2&\omega&\omega&\omega^2&1\\*[2mm]
1&1&\omega&\omega^2&\omega&\omega^2\\*[2mm]
1&\omega&1&\omega^2&\omega^2&\omega\\*[2mm]
1&\omega^2&\omega^2&1&\omega&\omega\end{array}\right]\label{A03}\end{eqnarray}
with the spectral function 
\begin{eqnarray}
f(A_{03})=(x^2-1)(x^4+\frac{1-\omega}{\sqrt{6}} x^3+x^2+\frac{\omega+2}{\sqrt{6}}x+1)
\label{sp03}\end{eqnarray}

From relations (\ref{sp01}) and   (\ref{sp03}) we observe that $f(A_{01})=f(A_{03})$, which means that from the six standard forms one gets only two different spectra given by equations (\ref{sp01}) and (\ref{sp02}). Thus at this level there are only two matrices that are not unitary equivalent.

The $A_3$ matrix has the form
\begin{eqnarray}
A_{3}(\omega)=\left[ \begin{array}{rrrrrr}
1&1&1&1&1&1\\*[2mm]
1&1&\omega&\omega^2&\omega&\omega^2\\*[2mm]
1&\omega&1&\omega^2&\omega^2&\omega\\*[2mm]
1&\omega^2&\omega^2 &1&\omega&\omega\\*[2mm]
1&\omega&\omega^2&\omega&1&\omega^2\\*[2mm]
1&\omega^2&\omega&\omega&\omega^2&1\end{array}\right]\label{Ag10}\end{eqnarray}

The choice $x = \omega^2,\,y = \omega,\,z = 1  $ leads to 
\begin{eqnarray}
A_{40}=\left[ \begin{array}{rrrrrr}
1&\omega^2&\omega&\omega &\omega^2&1 \\*[2mm]
\omega^2 &1&\omega&\omega^2&\omega&1 \\*[2mm]
\omega^2&\omega^2&\omega^2& \omega^2&\omega^2&\omega^2\\*[2mm]
1&\omega^2&1&\omega^2&\omega&\omega\\*[2mm]
\omega^2 &1&1&\omega &\omega^2&\omega \\*[2mm]
1&1&\omega^2 &\omega^2 &\omega&\omega^2\end{array}\right]\label{Ag9}\end{eqnarray}
whose  spectral function is the following
\begin{eqnarray}
 f(A_{40}) =  x^6+\sqrt{\frac{2}{3}}(1+2\omega)x^5-\frac{\omega+2}{2}x^4-\frac{1+2\omega}{\sqrt{6}} x^3+\frac{1-\omega}{2}x^2+\sqrt{\frac{2}{3}}(1+2\omega)x -1
\label{sp40}\end{eqnarray}

The choice $x = \omega^2,\,y = 1,\, z = \omega$ leads to $A_{50}$ matrix whose form is
\begin{eqnarray}
A_{50}=\left[ \begin{array}{rrrrrr}
\omega&\omega^2&1&1 &\omega^2&\omega \\*[2mm]
\omega^2 &\omega&1&\omega^2&1&\omega \\*[2mm]
\omega^2&\omega^2&\omega^2& \omega^2&\omega^2&\omega^2\\*[2mm]
\omega&\omega^2&\omega&\omega^2&1&1\\*[2mm]
\omega^2&\omega&\omega&1&\omega^2&1\\*[2mm]
\omega &\omega &\omega^2&1&1&\omega^2 
\end{array}\right]\label{Ag12}\end{eqnarray}
and its spectral function is given by
\begin{eqnarray}
f(A_{50})= x^6+\sqrt{\frac{2}{3}}(\omega+2)x^5+\frac{\omega+2}{2}x^4+\frac{2\omega+1}{\sqrt{6}} x^3+\frac{\omega-1}{2}x^2+\sqrt{\frac{2}{3}}(\omega-1)x -1
\label{sp50}\end{eqnarray}

The parameters $x = 1,\,y = \omega^2,\,z = \omega$ give rise to the matrix 
\begin{eqnarray}
A_{60}=\left[ \begin{array}{rrrrrr}
\omega&1&\omega^2&\omega^2&1 &\omega \\*[2mm]
1&\omega &\omega^2&1&\omega^2&\omega \\*[2mm]
1&1&1&1&1&1 \\*[2mm]
\omega&1&\omega&1&\omega^2&\omega^2\\*[2mm]
1&\omega&\omega&\omega^2&1&\omega^2\\*[2mm]
\omega &\omega &1&\omega^2&\omega^2 &1
\end{array}\right]\label{Ag14}\end{eqnarray}
with the spectral function 
\begin{eqnarray}
f(A_{60})= x^6-\sqrt{\frac{2}{3}}(\omega+2)x^5+\frac{2\omega+1}{2}x^4-\frac{2\omega+1}{\sqrt{6}} x^3+\frac{2\omega+1}{2}x^2-\sqrt{\frac{2}{3}}(\omega-1)x -1
\label{sp60}\end{eqnarray}

From the above examples one sees that $f(A_{i0})\, \ne\, f(A_{j0}) $, $\,\,i\,\ne\,j$, $\,\,i,\,j = 1,\dots,6$.

On the other hand all the  six matrices $A_i$ are different since $A_i\neq A_j$, $i\neq j$, (we have written only three of them), have the same spectrum given by relation (\ref{Sp(A1)}), thus they are {\em unitary equivalent}, and  the matrix (\ref{Ag2}) {\em is indeed an isolated one}. Why happens so\,?   An explanation could be the following:
the above  matrices $A_i$ can be transformed into symmetric matrices by the choice $\omega = a$, where {\em a} is a real number. For example $A_2$ has the new symmetric form
\begin{eqnarray}
A_{2}(a)=\left[ \begin{array}{rrrrrr}
1&1&1&1&1&1\\*[2mm]
1&1&a^2&a&a^2&a\\*[2mm]
1&a^2&1&a&a&a^2\\*[2mm]
1&a&a&1&a^2&a^2\\*[2mm]
1&a^2&a&a&1&a^2\\*[2mm]
1&a^2&a&a^2&a&1\end{array}\right]\end{eqnarray}
and its spectrum is given by

\begin{eqnarray}
Sp(A_{2}(a)) = \left[\frac{1+a+a^2\pm\sqrt{a^2(1+a^2)+5}}{\sqrt{6}}, \left(\frac{\pm\sqrt{5 a^2(a-1)^2}+2-a(1+a)}{2\sqrt{6}}\right)^2\right]\end{eqnarray}
and all the $A_i(a)$ matrices lead to the same spectrum, even if the matrices do not coincide. This thing could be interpreted as a higher symmetry of Agaian matrix since now $a$ is an arbitrary real number.

\section{Conclusion}

All the six forms of  $A_{i0}(\omega)$ matrices, $i\,= 1,\dots, 6$, obtained from matrix (\ref{Ag}) are {\em not unitary equivalent}. When they are  brought to their standard form only two of them are not unitary equivalent. 

The multiplication of $A_{i0}(\omega)$ matrices by diagonal matrices generated by their first row and column leads to a  higher symmetry.

Thus from the unitary equivalence 
 point of view we have {\em six} Agaian matrices of the form (\ref{Ag1}),  (\ref{Ag5}), etc, which are not unitary equivalent, and  only {\em two} Agaian matrices of the form (\ref{Ag3}), (\ref{Ag7}), which are brought to their standard form and which  are not unitary equivalent between themselves.

The final conclusion is that all the six form of the Agaian matrix whose first diagonal entries are equal to 1 are  equivalent, irrespective of the equivalence method which is used for. This shows that the Agaian matrix is isolated.
\begin{acknowledgments}
It a pleasure to thank Karol  $\dot{\rm Z}$yczkowski who suggested an analysis of the isolation problem for the Agaian matrix among the 6-dimensional matrices. 

We acknowledge a partial  support from Project PN09370102/2009.\end{acknowledgments}

\end{document}